\documentclass[11pt]{article}
\usepackage{amsmath, amssymb, amsthm}

\date{}
\title{\vspace{-0.8cm}Large induced trees in $K_r$-free graphs}
\author{
Jacob Fox \thanks{Department of Mathematics,
Princeton University, Princeton, NJ 08544. E-mail:
jacobfox@math.princeton.edu.
Research supported in part by an NSF Graduate Research Fellowship and a Princeton Centennial Fellowship}
\and
Po-Shen Loh \thanks{Department of Mathematics,
Princeton University, Princeton, NJ 08544. E-mail:
ploh@math.princeton.edu.
Research supported in part by a Fannie and John Hertz Foundation Fellowship, an NSF Graduate
Research Fellowship, and a Princeton Centennial Fellowship.}
\and
Benny Sudakov \thanks{Department of Mathematics, UCLA, Los Angeles, CA 90095. E-mail:
bsudakov@math.ucla.edu.
Research supported in part by NSF CAREER award DMS-0812005, and a USA-Israeli BSF grant.
}
}

\newtheorem{theorem}{Theorem}
\newtheorem{lemma}{Lemma}
\newtheorem{definition}{Definition}
\newtheorem*{claim*}{Claim}

\newcommand{\pr}{\mathbb{P}}
\newcommand{\bin}{\text{Bin}}

\begin{document}
\maketitle

\begin{abstract}
  For a graph $G$, let $t(G)$ denote the maximum number of vertices in
  an induced subgraph of $G$ that is a tree.  In this paper, we study
  the problem of bounding $t(G)$ for graphs which do not contain a
  complete graph $K_r$ on $r$ vertices.  This problem was posed twenty
  years ago by Erd\H{o}s, Saks, and S\'os.  Substantially improving
  earlier results of various researchers, we prove that every
  connected triangle-free graph on $n$ vertices contains an induced
  tree of order $\sqrt{n}$.  When $r \geq 4$, we also show that $t(G)
  \geq \frac{\log n}{4 \log r}$ for every connected $K_r$-free graph
  $G$ of order $n$.  Both of these bounds are tight up to small
  multiplicative constants, and the first one disproves a recent
  conjecture of Matou\v{s}ek and \v{S}\'amal.
\end{abstract}

\section{Introduction}

For a graph $G$, let $t(G)$ denote the maximum number of vertices in
an induced subgraph of $G$ that is a tree.  The problem of bounding
$t(G)$ in a connected graph $G$ was first introduced twenty years ago
by Erd\H{o}s, Saks, and S\'os \cite{ESS}.  Clearly, to get a non-trivial
result one must impose
some conditions on the graph $G$, because, for example, the complete
graph contains no induced tree with more than 2 vertices.  In their
paper, Erd\H{o}s, Saks, and S\'os studied the relationship between
$t(G)$ and several natural parameters of the graph $G$.
They were able to obtain asymptotically tight bounds on $t(G)$ when
either the number of edges or the independence number of $G$ were
known.

Erd\H{o}s, Saks, and S\'os also considered the problem of estimating
the size of the largest induced tree in graphs with no $K_r$ (complete
graph on $r$ vertices).  Let $t_r(n)$ be the minimum value of $t(G)$
over all connected $K_r$-free graphs $G$ on $n$ vertices.  In
particular, for triangle-free graphs, they proved that
\begin{displaymath}
  \Omega\left( \frac{\log n}{\log \log n} \right) \ \leq \ t_3(n) \ \leq \ O(\sqrt{n} \log n),
\end{displaymath}
and left as an interesting open problem the task of closing the wide
gap between these two bounds.

The first significant progress on this question was made only recently
by Matou\v{s}ek and \v{S}\'amal \cite{MS}, who actually came to the
problem of estimating $t_3(n)$ from a different direction.  Pultr had
been studying forbidden configurations in Priestley spaces \cite{BP},
and this led him to ask in \cite{P} how large $t(G)$ could be for
connected bipartite graphs $G$.  Let $t_B(n)$ be the minimum value of
$t(G)$ over all connected bipartite graphs on $n$ vertices.  It is
clear that $t_3(n) \leq t_B(n)$, so the result of Erd\H{o}s, Saks, and
S\'os immediately gives a lower bound on $t_B(n)$.

Motivated by Pultr's question, Matou\v{s}ek and \v{S}\'amal studied
$t_B(n)$ and $t_3(n)$.  They found the following nice construction
which shows that $t_3(n) \leq t_B(n) < 2\sqrt{n}+1$.  Let $m =
\sqrt{n}$, and consider the graph with parts
$V_{-m+1}$, $V_{-m+2}$, \ldots, $V_{m-1}$, where $|V_i| = m-|i|$, and
each consecutive pair of parts $(V_i, V_{i+1})$ induces a complete
bipartite graph.  This graph is clearly bipartite with $m^2 = n$
vertices, and it is easy to see that every induced tree in it has at
most $2m-1$ vertices.  On the other hand, Matou\v{s}ek and \v{S}\'amal
were able to improve the lower bound on $t_B(n)$ and $t_3(n)$, showing
that $t_3(n) \geq e^{c \sqrt{\log n}}$ for some constant $c$.
Furthermore, they also proved that if there was even a single value of
$n_0$ for which $t_3(n_0) < \sqrt{n_0}$, then in fact $t_3(n) \leq
O(n^\beta)$ for some constant $\beta$ strictly below $1/2$.  The above
fact led Matou\v{s}ek and \v{S}\'amal to conjecture that the true
asymptotic behavior of $t_3(n)$ was some positive power of $n$ which
is strictly smaller than $1/2$.

Our first main result essentially solves this problem.  It determines
that the order of growth of both $t_3(n)$ and $t_B(n)$ is precisely
$\Theta(\sqrt{n})$, disproving the conjecture of Matou\v{s}ek and
\v{S}\'amal.

\begin{theorem}
  \label{thm:r=3}
  Let $G$ be a connected triangle-free graph on $n$ vertices.  Then
  $t(G) \geq \sqrt{n}$.
\end{theorem}

Furthermore, our approach can also be used to give asymptotically
tight bounds on the size of the largest induced tree in $K_r$-free
graphs for all remaining values of $r$.  In their original paper,
Erd\H{o}s, Saks, and S\'os gave an elegant construction which shows
that $t_r(n)$ for $r \geq 4$ has only logarithmic growth.  Indeed, let
$T$ be a balanced $(r-1)$-regular tree, that is, a rooted tree in which
all non-leaf vertices have degree $r-1$ and the depth of any two
leaves differs by at most 1.  Then the line graph\footnote{The
  vertices of $L(T)$ are the edges of $T$, and two of them are
  adjacent if they share a vertex in $T$.}  $L(T)$ is clearly
$K_r$-free, and one can easily check that induced trees in $L(T)$
correspond to induced paths in $T$, which have only logarithmic
length.  Optimizing the choice of the parameters in this construction, one can show that
$t_r(n) \leq \frac{2 \log (n-1)}{\log (r-2)} + 2$.  On the other hand,
using Ramsey Theory, Erd\H{o}s, Saks, and S\'os also showed that
$t_r(n) \geq \frac{c_r \log n}{\log \log n}$, where $c_r$ is a
constant factor depending only on $r$.  Our second main result closes
the gap between these two bounds as well, and determines the order of
growth of $t_r(n)$ up to a small multiplicative constant.

\begin{theorem}
  \label{thm:r>3}
  Let $r \geq 4$, and let $G$ be a connected graph on $n$ vertices
  with no clique of size $r$.  Then $t(G) \geq \frac{\log n}{4 \log
    r}$.
\end{theorem}

One can also study induced {\em forests} rather than trees in $K_r$-free graphs.
Let $f_r(n)$  be the maximum number such that every $K_r$-free
graph on $n$ vertices contains an induced forest with at least
$f_r(n)$ vertices. Trivially we have $f_r(n) \geq t_r(n)$, but it appears that
the size of the maximum induced forest in a graph is more closely related to
another parameter. The {\em independence number}
$\alpha(G)$ of a graph is the size of the
largest independent set of vertices in $G$. 
Since an independent set is a forest and every
forest is bipartite,  the size of the largest
induced forest in a graph $G$ is at least $\alpha(G)$ and at most
$2\alpha(G)$. Using the best known upper bound for
off-diagonal Ramsey numbers \cite{AKS}, for fixed $r \geq 3$ and all
$n$ we have $f_r(n) \geq c n^{\frac{1}{r-1}}\log^{\frac{r-2}{r-1}} n$
for some positive constant $c$. Hence, $f_3(n)$ is larger than $t_3(n)$
by a factor of
$c \sqrt{\log n}$. Furthermore, for fixed $r>3$,
$f_r(n)$ and $t_r(n)$ behave very differently, as $f_r(n)$ is
polynomial in $n$ while $t_r(n)$ is only logarithmic. This
demonstrates that in $K_r$-free graphs the largest
guaranteed induced forest is much larger than the largest guaranteed
induced tree.

We close this introduction by mentioning some related research.  Our
work considers the Ramsey-type problem of finding either a clique or
a large induced tree.  The similar problem of finding an induced
copy of a \emph{particular}\/ tree $T$ in a $K_r$-free graph was
independently raised by Gy\'arf\'as \cite{G-conj} and Sumner
\cite{Sum}.  They conjectured that for any fixed integer $r$ and
tree $T$, any graph with sufficiently large chromatic number
(depending on $r$ and $T$) must contain either an $r$-clique or an
induced copy of $T$.  Note that the essential parameter for the
graph $G$ is now the chromatic number and not the number of
vertices.  Indeed, a complete bipartite graph has no clique of size
3, but contains only stars as induced subtrees.  This conjecture is
widely open, although some partial results were obtained in
\cite{G-path, GST, KP, Sc}.

Induced trees were also studied in the context of sparse random
graphs.  This line of research was started by Erd\H{o}s and Palka
\cite{EP}, who conjectured that for any constant $c>1$, the random
graph $G(n, c/n)$ would with high probability contain an induced tree
of order $\gamma(c)n$.  This was solved by Fernandez de la Vega
\cite{F1}, and other variants of this result were obtained in
\cite{KR, F2, FJ, LP, Sue}.  In another regime, when the edge
probability is $p = c \log n /n$, Palka and Ruci\'nski \cite{PR}
showed that the largest induced tree in $G(n, p)$ has size $\Theta(n
\log \log n / \log n)$ with high probability.

The rest of this paper is organized as follows.  In Section
\ref{sec:r=3} we discuss the proof of Theorem \ref{thm:r=3}, and show
how to reduce it to an abstract optimization problem on certain
bipartite graphs with weights on the vertices.  The solution of this
problem is provided in the following section.  In Section
\ref{sec:r>3}, we show how to extend our argument to the case of
$K_r$-free graphs with $r \geq 4$, and prove our second result,
Theorem \ref{thm:r>3}.  The final section of the paper contains a few
concluding remarks.  Throughout our paper, we will omit floor and
ceiling signs whenever they are not essential, to improve clarity of
presentation.

\section{Triangle-free graphs}
\label{sec:r=3}

The main idea in the proof of Theorem \ref{thm:r=3} is to use
induction to prove a slightly\footnote{We will discuss the relative
  strength of this statement in detail in our concluding remarks.}
stronger statement.  Instead of finding a single induced tree, we show
that no matter which vertex $v$ of the graph we choose, there exists a
large induced tree which contains $v$.  More precisely, we prove that
any connected, triangle-free graph with $n+1$ vertices contains an
induced tree of size $\sqrt{n}+1$ through any given vertex.

This is obviously true for $n=1$, which serves as the base of our
induction.  It remains to prove the statement for general $n \geq 2$,
while assuming its truth for all smaller values of $n$.  So, let $G =
(V,E)$ be an arbitrary connected triangle-free graph with $n+1$
vertices, and fix an arbitrary vertex $v \in V$.  We will find a large
induced tree through $v$.  Note that since $G$ is triangle-free,
$\{v\} \cup N(v)$ induces a star.  Therefore, we may assume that the
size of the neighborhood satisfies $|N(v)| < \sqrt{n}$, or else we are
done.

Consider the subgraph of $G$ induced by $V \setminus (\{v\} \cup N(v))$.
It decomposes into connected components, whose vertex sets we call $V_1,
\ldots, V_m$.  Now suppose that we could find a subset $U \subset
N(v)$, and a subset $I \subset [m]$, with the following properties:
\begin{description}
\item[(i)] For each $i \in I$, there is exactly one $u \in U$ which is
  adjacent to at least one vertex in $V_i$.  Let us denote this vertex by $u(i)$.
\item[(ii)] The sum $\sum_{i \in I} \sqrt{|V_i|}$ is at least
  $\sqrt{|V_1 \cup \ldots \cup V_m|}$.
\end{description}
Then, for each $i \in I$, we could apply the induction hypothesis to
the connected subgraph of $G$ induced by $\{u(i)\} \cup V_i$.  This
would give an induced tree $T_i$ containing $u(i)$, of size
$1+\sqrt{|V_i|}$.  Furthermore, it is easy to see that the union of
$\{v\}$ with all of the above constructed trees $T_i$ is also an
induced tree.  Indeed, since each $V_i$ is a maximal connected
component, there are no edges between the $T_i$, and since $G$ is
triangle-free, there are no edges inside $U \subset N(v)$.  Therefore,
we will have an induced tree with total size at least:
\begin{eqnarray*}
  |\{v\}| + |\{u(i) : i \in I\}| + \sum_{i \in I} \sqrt{|V_i|} &\geq& 1 + 1 + \sqrt{|V_1 \cup \ldots \cup V_m|} \\
  &=& 2 + \sqrt{|V \setminus (\{v\} \cup N(v)|)} \\
  &\geq& 2 + \sqrt{(n+1) - 1 - \sqrt{n}} \\
  &\geq& 1 + \sqrt{n},
\end{eqnarray*}
as desired.  Thus, the following abstract lemma completes the proof.

\begin{lemma}
  \label{lem:nonuniform-ineq}
  Consider a bipartite graph with sides $A$ and $B$, with the property
  that each vertex in $B$ has degree at least 1.  Let each vertex $i \in
  B$ have an associated weight $w_i \geq 0$.  We call a subset $H
  \subset A \cup B$ \textbf{admissible} if each vertex $v \in B \cap
  H$ has exactly one neighbor in $A \cap H$.  Then there exists an
  admissible $H$ with $\sum_{i \in B \cap H} \sqrt{w_i} \geq
  \sqrt{\sum_{i \in B} w_i}$.
\end{lemma}

The connection between this lemma and our required selection of $I
\subset [m]$ and $U \subset N(v)$ is clear.  The sides $A$ and $B$
correspond to the sets $N(v)$ and $[m]$, respectively, and the weights
$w_i$ are precisely the sizes of the connected components $|V_i|$.
The requirement that each vertex in $B$ has degree at least 1 is
satisfied by the fact that $G$ is connected, and so each component
$V_i$ has at least one neighbor in $N(v)$.  Therefore, this lemma will
indeed complete the proof of Theorem \ref{thm:r=3}.

\section{Main lemma}
\label{sec:main-lemma}

Before proving our main lemma, Lemma \ref{lem:nonuniform-ineq}, let us
discuss an easy special case which we will actually need later in our
study of $K_r$-free graphs when $r \geq 4$.  Observe that if the
weights $w_i$ in Lemma \ref{lem:nonuniform-ineq} were roughly equal,
then one way to control the objective $\sum_{i \in B \cap H}
\sqrt{w_i}$ would be to find a lower bound on $|B \cap H|$.  This
motivates the following claim, which we record for later use.

\begin{lemma}
  \label{lem:uniform-sqrt}
  Consider a bipartite graph with sides $A$ and $B$, with the property
  that each vertex in $B$ has degree at least 1.  We still call a
  subset $H \subset A \cup B$ \textbf{admissible} if each vertex $v
  \in B \cap H$ has exactly one neighbor in $A \cap H$.  Then there
  exists an admissible $H$ with $|B \cap H| \geq \sqrt{|B|}$.
\end{lemma}

\noindent \textbf{Proof.}\, The key observation, which we will also
use in the proof of Lemma \ref{lem:nonuniform-ineq}, is that we may
assume that every $v \in A$ has some vertex $w \in B$ which is
adjacent only to $v$.  Indeed, suppose this is not the case, and every
neighbor of $v$ has additional neighbors in $A$.  Then,
deleting $v$ from $A$ will not break the hypothesis of the lemma.
Therefore, after repeatedly performing this reduction, we obtain a
bipartite graph in which every vertex $v \in A$ has a neighbor in $B$
that sees only $v$.  Notice that this implies that there is an induced
matching between $A$ and some subset $B' \subset B$.

If $|A| \geq \sqrt{|B|}$, then the induced matching immediately yields
an admissible set $H = A \cup B'$ which satisfies the assertion.  On
the other hand, when $|A| < \sqrt{|B|}$, there is a vertex in $A$ with
degree at least $\sqrt{B}$.  Indeed, since every vertex in $B$ has
degree at least 1, the total number of edges in the bipartite graph is
at least $|B|$, and therefore some vertex $v \in A$ has degree $\geq
|B|/|A| > \sqrt{B}$.  The induced star $H = \{v\} \cup N(v)$ provides
the desired admissible set. \hfill $\Box$

\vspace{2mm}

We pause now to remark that Lemma \ref{lem:uniform-sqrt} is far from
being sharp.  In fact, it is always possible to find an admissible $H$
with $|H \cap B| \geq \Omega(|B|/\log |B|)$, and this is tight.
Although we do not need this result for our proof we sketch it here
for the sake of completeness.  By the reduction in the proof of Lemma
\ref{lem:uniform-sqrt}, we may assume that there is an induced
matching between $A$ and a subset $B' \subset B$.  In particular, this
implies that all degrees in $B$ are at most $|A| = |B'| \leq |B|$.
The set of possible degrees $\{1, 2, \ldots, |B|\}$ is covered by the
family of $\log_2 |B|$ dyadic intervals $I_k = [2^k, 2^{k+1}]$, so
there must be some $I_k$ with the property that at least $|B|/\log_2
|B|$ vertices of $B$ have degrees in $I_k$.  Sample a random subset
$A' \subset A$ by taking each vertex independently with probability $p
= 2^{-k-1}$, and let $B''$ be the set of all vertices in $B$ that are
adjacent to exactly one vertex in $A'$.  It is clear that $H = A' \cup
B''$ is admissible, so it remains to control $|B''|$.  Any vertex $v
\in B$ has probability exactly $\pr[\bin(d(v), 2^{-k-1}) = 1] =
d(v)2^{-k-1}(1-2^{-k-1})^{d(v)-1}$ of being chosen for $B''$.  Since
$I_k = [2^k, 2^{k+1}]$, when $d(v) \in I_k$ this probability is
bounded from below by an absolute constant (one can take $1/8$).
Hence the expected size of $B''$ is at least $\Omega(|I_k|) \geq
\Omega(|B|/\log |B|)$, which implies that there must exist some choice
of $A'$ and $B''$ that satisfy this bound.

The following construction shows that this bound is asymptotically
tight.  Choose integers $m = 2^k$, let $A = \mathbb{Z} / m\mathbb{Z}$,
and let $B = B_0 \cup \ldots \cup B_k$, where each $B_i = \{b_{i,1},
\ldots, b_{i,m}\}$.  Let each $b_{i,j}$ be adjacent to precisely $\{i,
i+1, \ldots, i+2^j-1\} \in A$, where we reduce everything modulo $m$.
This has $|B| = (k+1)m = \Theta(m \log m)$, but it is not too
difficult to verify that any admissible $H$ has $|H \cap B| < 2m$.

\subsection{Proof of Lemma \ref{lem:nonuniform-ineq}}

Unfortunately, Lemma \ref{lem:uniform-sqrt} is insufficient in general
for our application, because in our triangle-free graph, the sizes of
the connected components $V_i$ of $V \setminus (\{v\} \cup N(v))$ may
differ wildly.  For this, we need its weighted variant, which we prove
in this section.  The main trick in the proof is to vary the weights, which leads
us to study the following function.

\begin{definition}
  Let $G$ be a bipartite graph with vertex set $A \cup B$.  For
  notational convenience, let the vertices of $B$ be named $\{1, 2,
  \ldots, m\}$.  Then, we define
  \begin{displaymath}
    F_G(w_1, \ldots, w_m) = \max_{\text{\rm admissible}\ H \subset A \cup B} \sum_{i \in
      B \cap H} \sqrt{w_i},
  \end{displaymath}
  where we still say that a nonempty subset $H \subset A \cup B$ is
  \textbf{admissible} when every vertex in $B \cap H$ has exactly one
  neighbor in $A \cap H$.
\end{definition}

Lemma \ref{lem:nonuniform-ineq} is thus equivalent to the statement
that $F_G(\mathbf{w}) \geq \sqrt{\sum_{i=1}^m w_i}$ for any collection
of $w_i \geq 0$.  Since this inequality is homogeneous in the $w_i$,
from now on we will always assume that the weights have been
normalized to sum to 1.  It then suffices to show that
$F_G(\mathbf{w}) \geq 1$ for all $\mathbf{w}$ that satisfy the
constraints $w_i \geq 0$, $\sum w_i = 1$.  Observe that this domain is
now compact, and $F_G$ is a maximum of a finite collection of
continuous functions, hence continuous.  Therefore, $F_G$ attains its
infimum on this domain, which we will denote $\min_{\mathbf{w}} F_G$.

So, suppose for the sake of contradiction that we have some graph $G$
of minimum order for which $\min_{\mathbf{w}} F_G < 1$.  Graph $G$
must have an induced matching between $A$ and some subset $B' \subset
B$, since otherwise we can use the same reduction argument as in the
proof of Lemma \ref{lem:uniform-sqrt} to obtain a contradiction to the
minimality of $G$.  Let $(w_1,\ldots,w_m)$ be a minimizing assignment
for $F_G$, satisfying the constraints $w_i \geq 0$ and $\sum w_i = 1$.
Note that actually all $w_i$ must be strictly positive, or else we
could delete a vertex $i$ with $w_i=0$ to obtain a proper induced
subgraph $G'$ of $G$ and a weight assignment $\mathbf{w}'$ for which
$F_{G'}(\mathbf{w}')<1$, again contradicting the minimality of $G$.

We now exploit the fact that we have cast our problem in a
continuous setting.  Let us study the effect of performing the
following perturbation on the weights.

\begin{description}
\item[Stage 1.] For each $i \in B'$, let $w_i' = w_i - \epsilon
  \sqrt{w_i}$.  For $i \not \in B'$, let $w_i' = w_i$.
\item[Stage 2.] To compensate for the fact that $\sum w_i' = 1 -
  \epsilon \sum_{j \in B'} \sqrt{w_j} < 1$, renormalize by scaling up every
  weight by the same proportion.  That is, for all $i \in B$, let
  \begin{displaymath}
    w_i'' = \frac{w_i'}{1 - \epsilon \sum_{j \in B'} \sqrt{w_j}}.
  \end{displaymath}
\end{description}

Note that since all $w_i > 0$, for all sufficiently small $\epsilon$,
all new $w_i'$ are still positive.  This perturbation is chosen in the
particular way because for small $\epsilon$ and $i \in B'$,
$\sqrt{w_i'} = \sqrt{w_i} - \frac{\epsilon}{2} + o(\epsilon)$.  This
is because it is easy to check that for every $x > 0$,
\begin{displaymath}
  \lim_{\epsilon \rightarrow 0} \frac{\sqrt{x-\epsilon\sqrt{x}} - (\sqrt{x}-\epsilon/2)}{\epsilon} = 0.
\end{displaymath}
So, the effect on the square root of each weight $w_i$ with $i \in B'$
is roughly the same, no matter what the weight is.

Now recall that the function $F_G(w_1, \ldots, w_m)$ is defined as the
maximum over all admissible $H$ of $\sum_{i \in B \cap H} \sqrt{w_i}$.
Since all $w_i \geq 0$ by definition, this is equal to the maximum
over all \emph{maximal}\/ admissible $H$, where this maximality is
defined with respect to set inclusion.  For brevity, let $M = F_G(w_1,
\ldots, w_m)$ be that maximum
and let $H$ be any such maximal admissible selection.

Note that $H$ must
intersect $B'$ (the subset of $B$ that has an induced matching to
$A$).  This is because any maximal admissible $H$ contains at least
one vertex from $A$, and that vertex's partner in $B'$ can be added to
$H$ while preserving admissibility.  In particular, the sum $\sum_{i
  \in B \cap H} \sqrt{w_i}$ includes at least one downwardly perturbed
weight from $B'$.  Therefore,
\begin{displaymath}
  \sum_{i \in B \cap H} \sqrt{w_i'}
  \leq \left( \sum_{i \in B \cap H} \sqrt{w_i} \right) - \frac{\epsilon}{2} + o(\epsilon)
  \leq M - \frac{\epsilon}{2} + o(\epsilon).
\end{displaymath}

The renormalization that converts $w_i'$ into $w_i''$ is particularly
simple to analyze.  Using the previous inequality and the observation
that $\sum_{j \in B'} \sqrt{w_j} \leq M$ (because $B' \cup A$ is an
induced matching, hence admissible):
\begin{displaymath}
  \sum_{i \in B \cap H} \sqrt{w_i''}
  \leq \frac{M - \frac{\epsilon}{2} + o(\epsilon)}{\sqrt{1 - \epsilon \sum_{j \in B'} \sqrt{w_j}}}
  \leq \frac{M - \frac{\epsilon}{2} + o(\epsilon)}{\sqrt{1 - \epsilon M}}.
\end{displaymath}

The final bound is independent of $H$, so if it were strictly smaller
than $M$, we would have $F_G(w_1'', \ldots, w_m'') < M = F_G(w_i,
\ldots, w_m)$, contradicting the minimality of $(w_1, \ldots, w_m)$.
Therefore, we must have:
\begin{eqnarray*}
  \frac{M - \frac{\epsilon}{2} + o(\epsilon)}{\sqrt{1 - \epsilon M}} &\geq& M \\
  M - \frac{\epsilon}{2} + o(\epsilon) &\geq& M\sqrt{1 - \epsilon M} \\
  M - \frac{\epsilon}{2} + o(\epsilon) &\geq& M \left(1 - \frac{\epsilon M}{2} + o(\epsilon M)\right) \\
  - \frac{\epsilon}{2} + o(\epsilon) &\geq& -\frac{\epsilon M^2}{2} + o(\epsilon M^2) \\
  1 - o(1) &\leq& M^2.
\end{eqnarray*}
In the final inequality, we used the fact that $M$ is fixed, and
therefore $o(M^2) = o(1)$.  Sending $\epsilon$ to zero, we conclude
that $F_G(\mathbf{w}) = M \geq 1$.  This contradicts our assumption
that $F_G(\mathbf{w}) = \min_{\mathbf{w}} F_G < 1$, so our proof is
complete. \hfill $\Box$

\vspace{2mm}

\noindent \textbf{Remark.}\, The following example shows that the
assertion of Lemma \ref{lem:nonuniform-ineq} no longer holds for any
exponent $\alpha > 1/2$.  Indeed, consider the following bipartite
graph.  For some very large $t$, let $A = \{a_1, \ldots, a_t\}$, let
$B = \{b_0, \ldots, b_t\}$, and connect each $a_i$ to $b_0$ and $b_i$.
Let the weight of $b_0$ be $1-t^{-1}$, and the weights of all other
vertices in $B$ be $t^{-2}$, so the total weight is 1.  It is easy to
see that the only maximal admissible sets in this graph are either a
star containing $b_0$ and some other $b_i$, or the induced matching
between $A$ and $B \setminus \{b_0\}$.  Since $\alpha > 1/2$ and $t$
is sufficiently large, we have in the first case that
$(1-t^{-1})^\alpha + (t^{-2})^\alpha = 1 - \alpha t^{-1} + o(t^{-1}) <
1$.  On the other hand, for the second admissible set, we only have $t
\cdot (t^{-2})^\alpha = t^{1 - 2\alpha} < 1$.

\section{$K_r$-free graphs}
\label{sec:r>3}

This section is devoted to the proof of Theorem \ref{thm:r>3}.  The
induction approach we used in Section \ref{sec:r=3} easily extends to
the case of $K_r$-free graphs when $r \geq 4$, and in fact the
argument becomes even simpler.  We prove that for any $r \geq 4$,
every connected $K_r$-free graph $G=(V,E)$ with $n+1$ vertices
contains an induced tree of size $\frac{\log n}{4 \log r}+1$ through
any particular vertex.  Note that since the logarithm appears in both
the numerator and denominator, its base is irrelevant.  The statement
is clearly true for $n=1$, which starts our induction.

Now, consider any $n \geq 2$, and suppose that the statement holds for
all smaller values of $n$.  Let $v \in V$ be an arbitrary vertex.  We
will find an induced tree of size $\frac{\log n}{4 \log r}+1$
containing $v$.  Recall the well-known fact from Ramsey Theory (see,
e.g., chapter 6.1 of \cite{Ramsey}) that any graph with $a^b \geq
{a+b-2 \choose a-1}$ vertices contains either a clique of size $a$ or
an independent set of size $b$.  This implies that the degree of $v$
must be less than $r^{\log n/4 \log r} = n^{1/4}$, or else we would
already be done.  Indeed, since $G$ is $K_r$-free, the neighborhood of
$v$ would then contain an independent set of size $\frac{\log n}{4 \log
  r}$.  The vertices of this set together with $v$ form an induced
star of the desired size.  The same argument also shows that every $w
\in N(v)$ has less than $n^{1/4}$ neighbors in $V \setminus (\{v\}
\cup N(v))$.  Otherwise, by the above discussion, we could find an
independent set $I \subset V \setminus (\{v\} \cup N(v))$ of size
$\frac{\log n}{4 \log r}$, all of whose vertices are adjacent to $w$.
Then, $v$, $w$, and $I$ will form a large induced tree containing $v$.

Let $V_1$, \ldots, $V_m$ be the vertex sets of the connected
components of the subgraph of $G$ induced by $V \setminus (\{v\}
\cup N(v))$.  Since $G$ is connected, each $V_i$ is adjacent to some
vertex in $N(v)$.  As we explained above, each vertex in $N(v)$ is
adjacent to fewer than $n^{1/4}$ sets $V_i$, so in particular $m <
|N(v)| n^{1/4} < n^{1/2}$.  We claim that all components $V_i$ have
size at most $\frac{n}{r^4}$.  Indeed, suppose that some $|V_i|$
exceeds $\frac{n}{r^4}$.  Let $u$ be a vertex in $N(v)$ which is
adjacent to at least one vertex in $V_i$.  Applying the induction
hypothesis to $\{u\} \cup V_i$, we find an induced tree $T$ through
$u$ of size $\frac{\log (n/r^4)}{4 \log r} + 1 = \frac{\log n}{4
\log r}$.  Then $\{v\} \cup T$ gives an induced tree of the desired
size.

Next, we show that there are more than $r^2$ indices $i$ for which
$|V_i| \geq \frac{\sqrt{n}}{r^2}$.  Indeed, if this were not the case,
then the total number of vertices in $V$ would be less than:
\begin{displaymath}
  |\{v\} \cup N(v)| + \sum_{i=1}^m |V_i|
  < 1 + n^{1/4} + m\cdot \frac{\sqrt{n}}{r^2} + r^2 \cdot \frac{n}{r^4}
  < 1 + n^{1/4} + 2 \cdot\frac{n}{r^2}.
  \leq 1 + n^{1/4} + \frac{n}{8}.
\end{displaymath}
This is less than $n+1 = |V|$ for all $n \geq 2$, which is a
contradiction.

Let $B$ be the above set of indices for which $|V_i| \geq
\frac{\sqrt{n}}{r^2}$, and let $A = N(v)$.  Consider the auxiliary
bipartite graph with sides $A$ and $B$ obtained by connecting $u \in
A$ with $i \in B$ if $u$ is adjacent to at least one vertex in $V_i$.
Applying Lemma \ref{lem:uniform-sqrt}, we find subsets $A' \subset A$
and $B' \subset B$ with $|B'| \geq \sqrt{|B|} > r$ such that for each
$i \in B'$ the component $V_i$ is adjacent to exactly one vertex in
$A' \subset N(v)$, which we denote $u(i)$.  In fact, $|B'| \geq r+1$
since both $|B'|$ and $r$ are integers.  Apply the induction
hypothesis to each $\{u(i)\} \cup V_i$ to find an induced tree $T_i$
containing $u(i)$ of size at least $\frac{\log |V_i|}{4 \log r} + 1$.

If all $u(i)$ are distinct, then we can find $u(i) \neq u(j)$ which
are not adjacent in $G$, because this is a set of at least $r+1$
vertices in a $K_r$-free graph.  Then, $\{v, u(i), u(j)\} \cup T_i
\cup T_j$ is an induced tree.  On the other hand, if there is some
$u(i) = u(j)$, then $\{v, u(i)\} \cup T_i \cup T_j$ is an induced
tree.  In either case, we find an induced tree containing $v$ of size at least
\begin{displaymath}
  |\{v, u(i), u(j)\}| + \frac{\log |V_i|}{4 \log r} + \frac{\log |V_j|}{4 \log r}
  \geq 2 + 2 \cdot \frac{\log (\sqrt{n}/r^2)}{4 \log r}
  = 1 + \frac{\log n}{4 \log r},
\end{displaymath}
as desired.  This completes the proof.  \hfill $\Box$

\section{Concluding remarks}

In this paper, we obtain a lower bound on the size of the largest
induced-tree in a $K_r$-free graph, which is tight up to a small
multiplicative constant.  Moreover, our proof shows that we can find
a large tree through any particular vertex in the graph.  It turns out
that this seemingly stronger result is equivalent up to a constant
factor to the original problem of finding one large tree.

\begin{claim*}
  Let $T$ be an induced tree in a connected graph $G$, and let $v$ be an
  arbitrary vertex.  Then $G$ has an induced tree of size $1 +
  \frac{|T|}{2}$ which contains $v$.
\end{claim*}

\noindent \textbf{Proof.}\, If $v$ is already in $T$, then there is
nothing to prove.  Otherwise, let $P = (v_1, v_2, \ldots, v_m)$ be a
shortest path between $v$ and $T$, with $v_1 = v$ and $v_m \in T$.  By
minimality of $P$, there are no edges between $\{v_1, \ldots,
v_{m-2}\}$ and $T$.  Let $e_1, \ldots, e_k$ be the edges connecting
$v_{m-1}$ and $T$, and let $t_1, \ldots, t_k$ be their endpoints in
$T$.  Since $T$ is a tree, by deleting some edges, we can partition it
into subtrees $T_1$, \ldots, $T_k$, such that each $T_i$ contains
$t_i$.  Consider the auxiliary graph on $k$ vertices, in which
vertices $i$ and $j$ are adjacent if there is an edge of $G$ between
$T_i$ and $T_j$.  Note that this graph also forms a tree, and
therefore is bipartite.  Hence, we can find two disjoint subsets $I
\cup J = [k]$ such that the collection of $T_i$ with $i \in I$ has no
edges crossing between them, and similarly, the collection of $T_j$
with $j \in J$ also has no crossing edges.  Therefore, the union of
$\{v_1, \ldots, v_{m-1}\}$ with either one of these two collections
will form an induced tree.  Clearly, both of these trees contain $v$,
and their union covers $T$.  Thus, one of them has size at least $1 +
\frac{|T|}{2}$. \hfill $\Box$

\vspace{2mm}

We also wish to remark that for the problem of finding a large induced
tree through every vertex of a triangle-free graph, one can improve
the $2\sqrt{n}$ upper bound of Matou\v{s}ek and \v{S}\'amal.  Indeed,
consider the following triangle-free graph on $n$ vertices.  Let $m =
\sqrt{2n}$, and take the graph with parts $V_0$, \ldots, $V_{m-1}$,
where $V_0 = \{v\}$, every other $|V_i| = m-i$, and each consecutive
pair of parts $(V_i, V_{i+1})$ induces a complete bipartite graph.
This is a bipartite (hence also triangle-free) graph with $1 +
m(m-1)/2 = (1+o(1))n$ vertices, but one can easily check that any
induced tree containing $v$ has at most $m = \sqrt{2n}$ vertices.  In
particular, this shows that any approach which guarantees a large tree
through every vertex of the graph cannot match Matou\v{s}ek and
\v{S}\'amal's upper bound.

In light of this discussion, we do not have a clear conjecture as to
what is the right constant in front of $\sqrt{n}$ in the problem of
finding a maximum induced tree in a triangle-free graph on $n$
vertices.  Nevertheless, as the upper and lower bounds are now so
close, perhaps there is a hope to bridge this gap with other methods.

\end{document}